\documentclass[letterpaper, 10 pt, conference]{ieeeconf} 
\IEEEoverridecommandlockouts
\renewcommand{\baselinestretch}{0.94}

\usepackage{graphicx}
\usepackage{comment}
\usepackage{amsmath}
\usepackage{amsthm,amssymb}
\usepackage[ruled,vlined]{algorithm2e}
\usepackage{tikz}
\usepackage{mathtools}
\usepackage{standalone}
\usepackage{caption}
\usepackage{subcaption}
\usepackage{dsfont}
\usepackage[font=footnotesize,labelfont=bf]{caption}

\usepackage{multirow}
\usepackage{lettrine}
\usepackage[noadjust]{cite}

\usetikzlibrary{patterns}

\newcommand{\dsqbL}{[\kern-0.15em[}
\newcommand{\dsqbR}{]\kern-0.15em]}
\newcommand{\cycle}{\texttt{cycle}}

\usetikzlibrary{arrows}
\graphicspath{ {./fig/} }

\title{\LARGE \bf 
Joint Pricing and Rebalancing of Autonomous \\ Mobility-on-Demand Systems}

\author{Salom\'{o}n Wollenstein-Betech$^{1}$, Ioannis Ch. Paschalidis$^1$, and Christos G. Cassandras$^1$ 
\thanks{*This work was supported in part by NSF under grants ECCS-1509084, DMS-1664644, CNS-1645681, IIS-1914792, and CMMI-1454737, by AFOSR under grant FA9550-19-1-0158, by ARPA-E's NEXTCAR program under grant DEAR0000796, by the MathWorks, by the ONR under grant N00014-19-1-2571, and by the NIH under grant 1R01GM135930.
We thank D. Sverdlin-Lisker for proofreading this paper.}
\thanks{$^{1}$ The authors are with the
Division of Systems Engineering, Boston University, Brookline, MA 02446 USA {\tt\small \{salomonw, cgc, yannisp\}@bu.edu}}
}

\newcommand{\be}[1]{\begin{equation}\label{#1}}
\newcommand{\benon}{\begin{equation*}}  
\newcommand{\bemuln}[1]{\begin{multline}\label{#1}}
\newcommand{\bemul}{\begin{multline*}}
\newcommand{\bee}{\begin{eqnarray*}}
\newcommand{\eee}{\end{eqnarray*}}
\newcommand{\been}[1]{\begin{eqnarray}\label{#1}}
\newcommand{\eeen}{\end{eqnarray}}
\newcommand{\began}[1]{\begin{gather}\label{#1}}
\newcommand{\bega}{\begin{gather*}}
\newcommand{\bealn}[1]{\begin{align}\label{#1}}
\newcommand{\beal}{\begin{align*}}
\newcommand{\bealatn}[2]{\begin{alignat}{#1}\label{#2}}
\newcommand{\bealat}{\begin{alignat*}}
\newcommand{\bexalatn}[1]{\begin{xalignat}\label{#1}}
\newcommand{\bexalat}{\begin{xalignat*}}

\def\bc{{\mathbf c}}

\def\br{{\mathbf r}}

\def\bu{{\mathbf u}}
\def\bv{{\mathbf v}}

\def\texitem#1{\par\smallskip\noindent\hangindent 25pt
               \hbox to 25pt {\hss #1 ~}\ignorespaces}

\newcommand{\scrA}{\mathcal{A}}
\newcommand{\scrB}{\mathcal{B}}

\newcommand{\scrF}{\mathcal{F}}
\newcommand{\scrG}{\mathcal{G}}

\newcommand{\scrN}{\mathcal{N}}

\newcommand{\scrP}{\mathcal{P}}

\newcommand{\scrR}{\mathcal{R}}

\newcommand{\scrU}{\mathcal{U}}

\newcommand{\blambda}{\boldsymbol{\lambda}}

\newcommand{\bmu}{\boldsymbol{\mu}}

\newtheorem{theorem}{Theorem}
\newtheorem{corollary}{Corollary}
\newtheorem{lemma}{Lemma}
\newtheorem{proposition}{Proposition}
\newtheorem{definition}{Definition}

\begin{document}
    \maketitle
    \begin{abstract}
This paper studies optimal pricing and rebalancing policies for Autonomous Mobility-on-Demand (AMoD) systems. We take a macroscopic planning perspective to tackle a profit maximization problem while ensuring that the system is load-balanced. We begin by describing the system using a dynamic fluid model to show the existence and stability of an equilibrium (i.e., load balance) through pricing policies. We then develop an optimization framework that allows us to find optimal policies in terms of pricing and rebalancing. We first maximize profit by only using pricing policies, then incorporate rebalancing, and finally we consider whether the solution is found sequentially or jointly. We apply each approach on a data-driven case study using real taxi data from New York City. Depending on which benchmarking solution we use, the joint problem (i.e., pricing and rebalancing) increases profits by 7\% to 40\%. 
\end{abstract}


    \section{Introduction} \label{sec:intro}

\lettrine{W}{ith} the rise of Mobility-on-Demand (MoD) services (e.g. Uber, Lyft, DiDi) and the rapid technological evolution of self-driving vehicles, we are closer to having Autonomous Mobility-on-Demand (AMoD) systems. A crucial step in the proper functioning of such a service is to define pricing, rebalancing and routing policies for the fleet of vehicles. This paper focuses on the first two issues, while the interested reader is directed to  \cite{wollensteinbetech2020congestionaware} for a discussion on routing and rebalancing.

\begin{figure}[t!]
        \centering
        \includegraphics[trim={0cm 0.5cm 0cm 4cm},clip, width=0.85\linewidth]{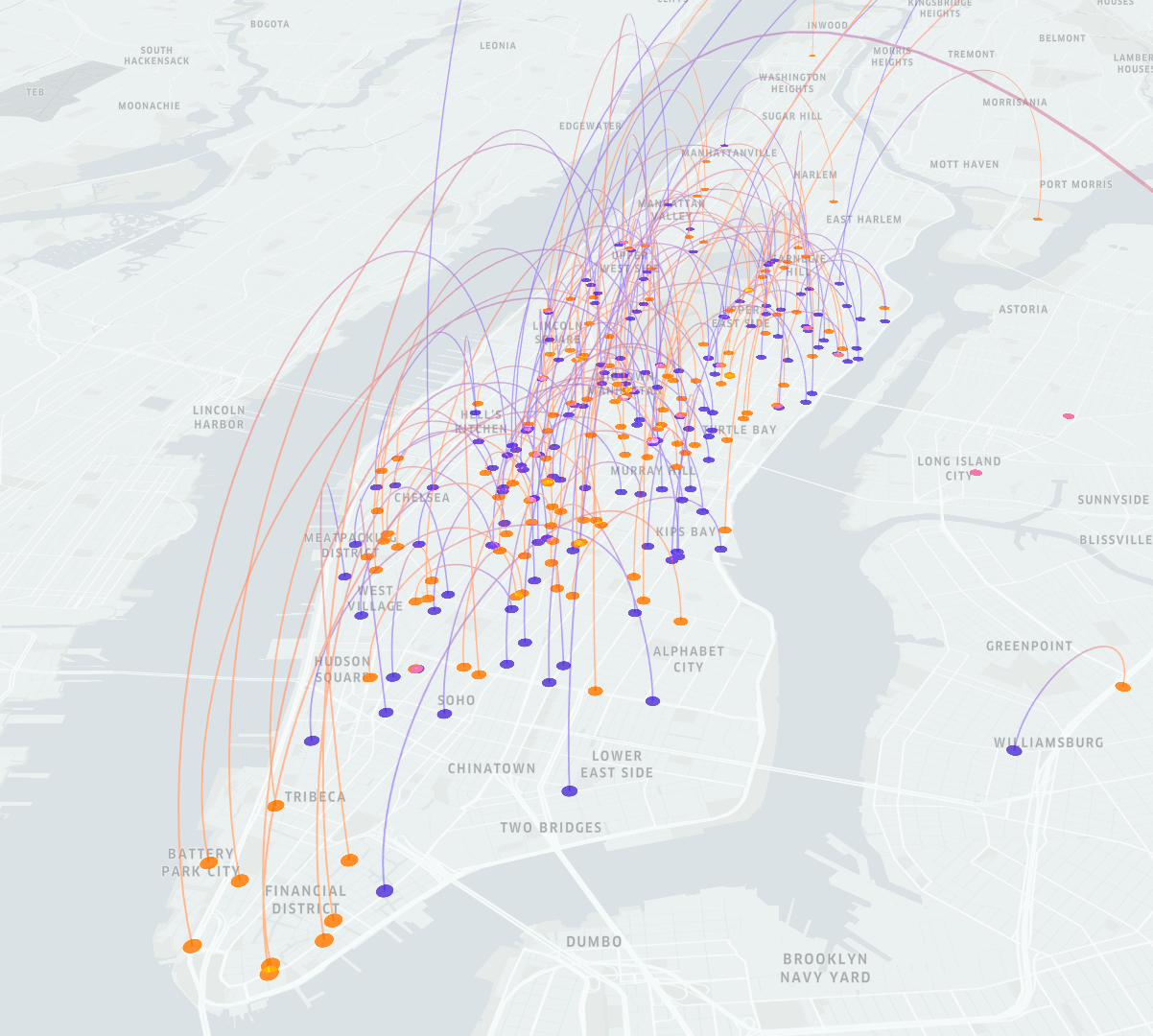}
        \caption{Requested taxi trips on January 15, 2015 10:38 a.m. in NYC. Blue and Orange circles represent origins and destinations respectively. One can observe that at this time, the Financial District (south) is an attractive destination but not origin. Hence, we expect taxis to rebalance to more attractive pickup locations.}
        \label{fig:amod}
        \vspace{-0.5cm}
    \end{figure}

Pricing policies play an important role as they modulate the inflow of customers traveling between regions in the network. As a result, the controller has the ability to choose prices such that the induced demand ensures a balanced load of customers and vehicles arriving at each location. 
Additionally, selecting prices enables the operator to modify demand such that the system can operate with smaller or larger fleet sizes. 
If we restrict a pricing policy to one that requires balancing the load in every node, we expect the solution to concentrate on balancing the network rather than choosing a set of prices to maximize profit. To give the pricing policy more flexibility, AMoD systems can leverage rebalancing policies, i.e., send empty vehicles from regions with excess supply of  vehicles to regions with excess demand of trips (see Fig.\ref{fig:amod}) with the objective of achieving higher profit.

\emph{Related Literature: }
Researchers have tackled the pricing problem using two main settings: one-sided, or two-sided markets depending on whether the MoD controller has full or limited control over the supply. In particular, one-sided markets assume full control over the vehicles \cite{banerjee2015pricing, turan2019dynamic}, whereas two-sided markets consider self-interested suppliers \cite{bimpikis2019spatial, banerjee2015pricing}. 
To the best of our knowledge, all these optimal pricing policies, except \cite{turan2019dynamic}, do not rebalance externally. Rather, they incentivize the supply (human drivers) to reallocate by the use of compensations. Our model differs from \cite{turan2019dynamic}, which uses a microscopic model and Reinforcement Learning techniques, by the level of abstraction performed. Alternatively to a microscopic model we employ a macroscopic (planning) model to assess the benefits of \emph{jointly} solving the pricing and rebalancing problem over other approaches.

The rebalancing literature has tackled the problem without the help of pricing incentives. For AMoD systems this problem has been studied using simulation models~\cite{swaszek2019load,HoerlRuchEtAl2018,LevinKockelmanEtAl2017}, queuing-theoretical models~\cite{ZhangPavone2016,IglesiasRossiEtAl2016}, network-flow models~\cite{PavoneSmithEtAl2012,RossiZhangEtAl2017} and it has also been studied jointly with routing schemes~\cite{SalazarTsaoEtAl2019,wollensteinbetech2020congestionaware}. In~\cite{swaszek2019load}, the rebalancing of an AMoD system is addressed using a data-driven parametric controller which is available for real-time implementation. Alternatively, in~\cite{PavoneSmithEtAl2012}, the rebalancing problem is studied using a steady-state fluid model which serves as a basis for this paper.

\emph{Key contributions: }
In this work we provide a theoretical framework to design optimal pricing policies for an AMoD provider. We analyze the system in the spirit of~\cite{pavone2012}, 
converting the problem into profit maximization rather than an operational cost minimization. Different from the existing methods in pricing, we consider the destination of a customer when designing the pricing policy. This allows the fleet controller to modulate demand in such a way that the system is balanced by solely adjusting prices. Additionally, we incorporate the rebalancing policy optimization framework in~\cite{pavone2012} and formulate a \emph{joint} optimization model. We compare this joint strategy with four different methodologies. First by only finding optimal prices, second by only rebalancing the fleet, third by sequentially solving the rebalancing and then pricing of the system, and fourth by jointly estimating pricing and rebalancing with a unique surge price by origin. We apply each approach to two case studies; one, with simulated data; and another, with real taxi data from New York City.

\emph{Organization:}
The paper is organized as follows. In Section \ref{sec:model} we introduce the fluid model consisting of queues of customers and vehicles at every region. In Section \ref{sec:well-poss-eq-stability}, we show that the system is well-posed and establish the existence of a load balance equilibrium through the selection of prices. We also obtain local stability results. Next, in Section \ref{sec:optimal-pricing}, we state the problems of optimal pricing, optimal rebalancing and the \emph{joint} formulation of these two. Then, we present case studies to assess the performance of the \emph{joint} formulation in Section \ref{sec:Experiments}. Finally, in Section \ref{sec:conclusion} we conclude.
    \section{Model} \label{sec:model}
In this section we present a steady-state deterministic fluid model to find optimal prices in an AMoD system while ensuring service to customers.
This  model is intended to serve as a relaxation of the corresponding stochastic queueing model where customers arrive according to a Poisson process and travel times are non-deterministic (usually exponentially distributed). The reason for making this relaxation is the flexibility it provides to perform analysis of the system. 

Consider a fully-connected network $\scrG=(\scrN, \scrA)$ where $\scrN$ is the set of nodes (regions)  $\scrN=\{1,...,N\}$  and $\scrA=\{(i,j) : i,j \in \scrN \times \scrN \}$ is the set of arcs. A customer requests a ride in region $i$, receives a transportation service from the AMoD platform, and is charged a price composed of the product of a \emph{base} and a  \emph{surge} price. The total price is  $p_{ij} = p^0_{ij} u_{ij}$ where $p^0_{ij}$, $u_{ij}$ are the base and surge prices, respectively, for traveling from node $i$ to $j$. Throughout the paper, we will use the surge price $u_{ij}$ as our control variable,  and we assume that $u_{ij}\geq 1$  as the platform is not willing to charge less than its base price.

We further assume that customers' arrival rate is a function of the surge price, namely $\lambda_{ij}(u_{ij}): \mathbb{R}_{\geq0} \mapsto \mathbb{R}_{\geq0}$ for a customer travelling from $i$ to $j$. This function is known as the \emph{willigness-to-pay} or the \emph{demand} function. Let the \emph{base demand} be  $\lambda^0_{ij}=\lambda_{ij}(1)$, i.e., the demand rate of customers when the surge price is at its minimum.

\begin{figure}[t]
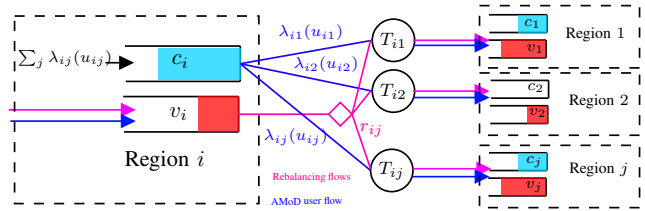

    \centering
    \includestandalone[trim={0 0 0 0},width=\linewidth]{systemDiagram}
    \caption{Customer traveling from $i$ to $j$ arrive to region $i$ at rate $\lambda_{ij}(u_{ij})$ and it takes $T_{ij}$ units of time to reach $j$. The AMoD provider plans a pricing policy $\bu$ and a rebalancing policy of empty vehicles $r_{ij}$ to serve their customers such that its profit is maximized. Note this is a fluid model as opposed to a discrete event system.}  
    \label{fig:system-diagram}
    \vspace{-0.5cm}
\end{figure}

As in~\cite{pavone2012}, we use a queueing model for this system with two queues per region. We let $c_i(t) \in \mathbb{R}_{\geq0}$ be the number of customers at region $i$ waiting to be assigned to a vehicle; and denote with $v_i(t) \in \mathbb{R}_{\geq0}$ the number of available vehicles waiting in region $i$ at time $t$. 
Moreover, the AMoD provider assigns vehicles to customers located in the same region at a service rate $\mu_i$. We assume that $\mu_i>\sum_j \lambda^0_{ij}$, meaning that the platform assigns vehicles to customers faster than the rate at which customers arrive. This assumption is required to avoid building large customer queues. For the purpose of this paper, we consider the rate vectors $\blambda=(\lambda_{ij}; \ \forall i,j \in \scrN )$ and $\bmu=(\mu_i; \ \forall i \in \scrN)$ to be invariant (we use bold notation to represent a vector containing all the variables sharing the same symbol). This allows us to analyze the steady-state solution of the system. Finally, we let $T_{ij}\in \mathbb{R}_{\geq0}$ be the travel time for a passenger to go from $i$ to $j$, which we assume to be fixed and not dependent on the routing decisions of the AMoD system (see Fig. \ref{fig:system-diagram}). 
To continue with our analysis, we make the following assumptions:\\

\textbf{Assumption 1.} The function $\lambda_{ij}(\cdot)$ is monotonically decreasing $\forall i,j \in \scrN$, i.e., as price increases, the demand rate decreases. 

\textbf{Assumption 2.} There exists a surge price $u^{\max}_{ij}$ for which $\lambda_{ij}(u^{\max}_{ij})=0, \quad \forall i,j \in \scrN$.

\subsection{Customer Dynamics}
\label{subsec:customer-dynamics}
Consider a customer queue $c_i$ for each region $i\in\scrN$ in the network. The queue dynamics are:
\begin{alignat}{3} \nonumber
    \dot{c_i} = & 
    \begin{dcases} 
        \sum\nolimits_j \lambda_{ij}(u_{ij}), & \text{if } v_i=0, \\ 
        0, & \text{if } v_i \geq 0 \text{ and } c_i = 0, \\
        \sum\nolimits_j \lambda_{ij}(u_{ij}) - \mu_{i} , &  \text{if } v_i \geq 0 \text{ and } c_i \geq 0.
    \end{dcases}
\end{alignat}
In order to express the customer dynamics with shorter notation we let $H(x)=\mathds{1}_{x>0}$ be an indicator function for positive values of $x$, and we use the following shorthand notation:
\begin{align*}
    \lambda_{ij} := \lambda_{ij}(u_{ij}), \quad
    \lambda_i := \sum\nolimits_j \lambda_{ij}, \quad
    v_i := v_i(t), \quad \\
    c_i := c_i(t), \quad
    v^i_j := v_j(t - T_{ji}), \quad
    c^j_i := c_j(t - T_{ji})
\end{align*}
where $\lambda_i$ is the total endogenous outgoing flow from node $i$; and $c^i_j$, $v^i_j$ are the rates of customer and vehicle arrivals at $t$ to region $i$ coming from $j$, respectively, i.e., $v^i_j$ is the rate of vehicles who departed from $j$ destined for $i$,  $T_{ji}$ time units prior to the current time $t$.
Then, we rewrite the customer dynamics in compact form can be written as follows:
\begin{align} \label{eq:customer-dynamics}
    \dot{c_i} = \lambda_{i} (1-H(v_i)) + (\lambda_{i} - \mu_i) H(c_i)H(v_i).
\end{align}
Note that as a result of using a fluid model, the variables denoting the number of customers in a region are real-valued. 

\subsection{Vehicle Dynamics} \label{subsec:vehicle-dynamics}
The outflow rate corresponding to vehicles departing station $i$ is given by:
\begin{alignat*}{3}
    \dot{v}_i^- = & 
        \begin{dcases}
            -\lambda_i, & \text{if } v_i \geq 0 \text{ and } c_i = 0, \\
            0, &  \text{if } v_i=0, \\
            -\mu_{i} , &  \text{if } v_i \geq 0 \text{ and } c_i \geq 0.
        \end{dcases}
\end{alignat*}
and more succinctly as 
\begin{equation} \label{eq:veh-arrival-rate}
    \dot{v}_i^- = - \lambda_i H(v_i) + (\lambda_i - \mu_i) H(v_i)H(c_i).
\end{equation}
In addition, the rate at which customer-carrying vehicles arrive at station $i$ is given by:
\begin{equation} \label{eq:veh-departure-rate}
    \dot{v}_i^+ = 
    \sum_j (\lambda_{ji}H(v^i_j) -
    (\lambda_{ji}-\mu_j)H(v^i_j)H(c^i_j) ).
\end{equation}
Therefore, we can write the vehicle dynamics in compact form by adding \eqref{eq:veh-arrival-rate} with \eqref{eq:veh-departure-rate} as
\begin{align} \label{eq:vehicle-dynamics}
    \dot{v_i}  =  - \lambda_{i} H(v_i) + (\lambda_{i} - \mu_i) H(c_i)H(v_i) &  \\
     +  \sum\limits_{j} (\lambda_{ji}H(v^i_j) - (\lambda_{ji} - \mu_j) &  H(c^i_j)H(v^i_j)) \nonumber.
\end{align}

Hence the global system dynamics are expressed by the following differential equations
\begin{subequations} \label{eq:system-dynamics}
\begin{flalign} 
    & \hspace{-0.8cm} \dot{c_i} = \lambda_{i} (1-H(v_i)) + (\lambda_{i} - \mu_i) H(c_i)H(v_i), \label{eq:system-dynamics-customers}\\
    & \hspace{-0.8cm} \dot{v_i}  =  - \lambda_{i} H(v_i) + (\lambda_{i} - \mu_i) H(c_i)H(v_i)  \label{eq:system-dynamics-veh}\\
    & \quad +  \sum\limits_{j} (\lambda_{ji}H(v^i_j) - (\lambda_{ji} - \mu_j)  H(c^i_j)H(v^i_j)) \nonumber.
\end{flalign}
\end{subequations}
which describe a non-linear, time-delayed, time-invariant, right-hand discontinuous system.  
    \section{Well posedness, Equilibrium and Stability} \label{sec:well-poss-eq-stability}

Similar to \cite{pavone2012}, we say that the system \eqref{eq:system-dynamics} is \emph{well posed} if two conditions are satisfied: \textit{(i)} for any initial condition, there exists a solution of the differential equations in \eqref{eq:system-dynamics}, and \textit{(ii)}, the number of vehicles in the system remain invariant over time. In order to analyze the model, we use the framework of Filippov solutions \cite{filippov2013differential}.
Let us now give a proposition for the well-posedness of the system:

\begin{proposition}[Well-posedness of fluid model] \label{prop:well-posedness}
    \textit{ 
    \begin{enumerate} 
        \item For every initial condition in the fluid model represented in \eqref{eq:system-dynamics}, there exist continuous functions $c_i(t): \mathbb{R}_{\geq 0} \mapsto \mathbb{R}_{\geq 0}$ and $v_i(t):\mathbb{R}_{\geq 0} \mapsto \mathbb{R}_{\geq 0}, \forall i \in \scrN$, satisfying the system of equations in the Fillipov sense. 
        \item For all $t>0$, the total number of vehicles is invariant and equal to $m=\sum_{i\in \scrN} v_i(0)$.
    \end{enumerate}
    }
\end{proposition}
\begin{pf}
For the first claim, we use the framework developed by \cite{haddad1981monotone}. In particular, we check that all assumptions and conditions of \cite[Thm II-1]{haddad1981monotone} are satisfied. This theorem, ensures the existence of Fillipov solutions to the time-delayed differential equations with discontinuous right-hand sides.
    
To prove the second claim, we study the dynamics of the vehicles in the system, which we separate into two categories: vehicles who are in transit $v_{ij}(t)$, and vehicles at a specific region $v_i(t)$. For the vehicles queued at $i$ we know their dynamics are as in \eqref{eq:system-dynamics-veh}. For the vehicles in transit, we let the total be
\begin{align*}
    v_{ij}(t)=\hspace{-3mm}\int\limits_{t-Tij}^t \hspace{-2mm} \lambda_{ij}H(v_i(\tau)) + (\lambda_{ij} - \mu_i) H(c_i(\tau))H(v_i(\tau)) \ d\tau,
\end{align*}
and their dynamics are
\begin{align*}
    \dot{v}_{ij}(t) =  \lambda_{ij}H(v_i) + (\lambda_{ij} - \mu_i) H(c_i)H(v_i) & \\
    - (\lambda_{ij}H(v^j_i) + (\lambda_{ij} - \mu_i)  &H(c^j_i)H(v^j_i)). 
\end{align*}
Moreover, the total number of vehicles in the system is $m(t)  = \sum_i v_i(t) + \sum_{ij} v_{ij}(t)$, with dynamics 
\begin{subequations}\label{eq:fleet-size-dynamics}
    \begin{align} 
        \dot{m}(t)& =   \sum\nolimits_i \dot{v}_i(t)+  \sum\nolimits_{ij} \dot{v}_{ij}(t)   \label{eq:fleet-size-dynamics-0},\\
        & = \sum\nolimits_{i} \big(-\lambda_{i} H(v_i) + (\lambda_{i} - \mu_i)  H(c_i)H(v_i)  \label{eq:fleet-size-dynamics-1} \\ 
        & \hspace{-1.4em} + \sum\nolimits_{j} \lambda_{ji}H(v^i_j) - (\lambda_{ji} - \mu_j) H(c^i_j)H(v^i_j)\big) + \sum\nolimits_{ij} \dot{v}_{ij},  \notag \\
        &=   \sum\nolimits_{ij} -\lambda_{ij} H(v_i) + (\lambda_{ij} - \mu_i) H(c_i)H(v_i)  \label{eq:fleet-size-dynamics-2} \\
        & \hspace{-1.4em} + \sum\nolimits_{ij} \lambda_{ji}H(v^i_j) - (\lambda_{ji} - \mu_j) H(c^i_j)H(v^i_j) + \sum\nolimits_{ij} \dot{v}_{ij}, \notag \\
        & =  0. 
    \end{align}
\end{subequations}
Note that to obtain the above result we have expanded the first sum term in \eqref{eq:fleet-size-dynamics-0} using \eqref{eq:system-dynamics-veh}, rearranged terms and found that $-\sum_i \dot{v}_i(t) =  \sum_{ij} \dot{v}_{ij}(t) \implies \dot{m}=0$, which implies that the fleet size remains invariant over time. 
\qed{}
\end{pf}

\subsection{Equilibria}
We say that the system is in equilibrium if customer queues (and therefore, waiting times) do not grow infinite. We show the existence of an equilibrium in the fluid model \eqref{eq:system-dynamics} when we control the prices of every origin-destination pair. Additionally, we show that by having the ability to control the prices, one can have find multiple equilibria for a desired fleet size, giving the flexibility to AMoD managers to operate the system at different levels.

\begin{theorem}[Existence of equilibria] \label{thm:existance-eq} 
    Let $\scrU$ be a set of prices $\bu$, such that when $\bu \in \scrU$:
    \begin{equation} \label{eq:equilibria}
        \sum\limits_{j} \lambda_{ij}(u_{ij}) - \lambda_{ji}(u_{ji}) = 0, \quad \forall i \in \scrN,
    \end{equation}
    and let 
    \begin{equation} \label{eq:min-vehicles}
        m_{\bu} := \sum\limits_{ij} T_{ij}\lambda_{ij}(u_{ij}).
    \end{equation}
   Then, if $\bu \in \scrU$, and $m > m_\bu$, an equilibrium exists with $\bc = 0$ and $\bv > 0$. Otherwise no equilibrium exists.
\end{theorem}

\begin{pf}
    Set $\dot{c}_i = 0 $ and $\dot{v}_i = 0 $ for all $i \in \scrN$. Then by using the customer system dynamics in \eqref{eq:system-dynamics-customers}, we have:
    \begin{equation} \label{eq:equlibria-dynamics-customers}
        \lambda_{i} = \lambda_{i}H(v_i) - (\lambda_{i} - \mu_i) H(c_i)H(v_i),
    \end{equation}
    and since $\lambda_i < \mu_i$, the above equation just has a solution if $c_i=0$ and $v_i>0$ for all $i \in \scrN$. Setting $\dot{v}_i = 0 $, and using the vehicle dynamics in \eqref{eq:system-dynamics-veh} we have 
    \begin{align} \label{eq:equlibria-dynamics-veh}
        0 = - & \lambda_{i} H(v_i) + (\lambda_{i} - \mu_i) H(c_i)H(v_i) \notag \\
        & + \sum\limits_{j} \lambda_{ji}H(v^i_j) - (\lambda_{ji} - \mu_j) H(c^i_j)H(v^i_j),
    \end{align}
    which combined with equation \eqref{eq:equlibria-dynamics-customers} and the fact that $\bc=0$ implies that
    \begin{equation} \label{eq:equilibria-Hv}
        0 = - \lambda_{i}  + \sum\limits_{j} \lambda_{ji}H(v_j).
    \end{equation}
    To arrive at \eqref{eq:equilibria-Hv}, we used the fact that in the stationary equilibrium $v_i$ and $c_i$ are constants and hence, there is no dependence on $t-T_{ij}$.
    
    Recall that for every equilibrium solution, we require $\bv > 0$ and thus $H(v_i)=1, \ \forall \in \scrN$. Therefore, a necessary condition for the existence of equilibria is that the prices $\bu$ are chosen such that
    \begin{align*}
        0 = - \lambda_{i}  + \sum\limits_{j} \lambda_{ji}, \quad \forall i \in \scrN,
    \end{align*}
    which proves the first statement.
    
    We now want to verify that the fleet size is large enough to maintain an equilibrium flow. Recall the fleet size dynamics $\dot{m}(t)$ when $\bc=0$ and $\bv>0$ in  \eqref{eq:fleet-size-dynamics}. Observe that to satisfy $\bv>0$, one needs to have a fleet size of at least $\sum\limits_{ij} T_{ij}\lambda_{ij}(u_{ij})$ vehicles which is the definition of $m_{\bu}$. This, mixed with its invariant property ($\dot{m}=0$), proves the claim. Conversely, if $m < m_{\bu}$ no equilibrium exists.
    \qed{}
    \end{pf}

\begin{lemma}[Existence of an equilibrium]
    The set $\scrU$ is never empty, hence, at least one equilibrium exists.
\end{lemma}
\begin{pf}
     We use the fact that there exists a price $u^{\max}_{ij}$ for which $\lambda_{ij}(u^{\max}_{ij}) = 0$ for all $i,j \in \scrN$.  Then, setting $\bu = \bu^{\max}$, implies that an equilibrium exists. This strategy means that we are not providing service to any request, nevertheless the equilibrium exists as we are unable to have an invariant fleet size. 
     \qed{}
\end{pf}

\begin{lemma}[Infinite number of equilibria]
    If there is a positive demand tour in the graph, then there exists an infinite number of price vectors $\bu$ which can steer the system to an equilibrium point.
\end{lemma}
\begin{pf}
    Assume that there exists at least one Eulerian tour (or \cycle) in the graph for which $\lambda^0_{ij} > 0$ for all $ (i,j) \in \cycle$. Then, let $\blambda^{\cycle}=\{\lambda^0_{ij}\ |\ (i,j) \in \cycle\}$ and the minimum rate on that tour be $\lambda^{\cycle}_{\min} = \min \{\lambda_{ij}\}_{(i,j) \in \cycle}$. Then by setting $u_{ij}=0$ for all $(i,j) \not \in \cycle$,  we can express the equilibrium condition as 
    \begin{equation}
        \sum\limits_{j:(i,j) \in \cycle} \hspace{-1.5em} \lambda_{ij}(u_{ij}) - \lambda_{ji}(u_{ji}) = 0, \ \ \forall i:(i,j) \in \cycle.
    \end{equation}
    Now, we use the fact that $\lambda_{ij}(u_{ij})$ is a monotonically decreasing function and we focus on $(i,j) \in \cycle$. Hence for all $\lambda_{ij}(u_{ij})>\lambda^{\cycle}_{\min}$ we can find a $u_{ij}$ such that $\lambda_{ij}(u_{ij}) = \lambda^{\cycle}_{\min}$. Then, extending this for higher prices on  $\lambda^{\cycle}_{\min}$ and using the same argument as before, we show that there exists a pricing strategy $\bu$ for which we can obtain an equilibrium with a tour demand rate with any value in the range $(0 ,\lambda^{\cycle}_{\min})$.
    \qed{}
\end{pf}
\\

These two lemmata imply that by incorporating an origin-destination pricing strategy, we can operate a mobility-on-demand service at equilibrium for any demand rate and with any fleet size. 

\begin{corollary}[Minimum number of vehicles in equilibria]
    The minimum number of vehicles to operate in an equilibrium induced by policy $\bu$ is at least
    \begin{align*}
        m > \underline{m} := \min_{\bu } m_{\bu}
    \end{align*}
    where $m_{\bu} := \sum\limits_{ij}T_{ij}\lambda_{ij}(u_{ij})$.
\end{corollary}
\begin{pf}
    This result follows directly from the last argument in the proof of Theorem \ref{thm:existance-eq}. \qed{}
\end{pf}

\subsection{Stability} \label{subsec:stability}

In this section we study local stability of the equilibria presented in the previous subsection. As an example, we look at cases when a disruptive change happens to the system, either because of an increase in customers or a decrease in the availability of vehicles.  
Let $\bu\in\scrU$ and assume $m_{\bu}>\underline{m}$. Then, we define the set of equilibria as 
\begin{flalign}
    & \Upsilon_{\bu} := \{(\bc, \bv) \in \mathbb{R}^{2N} \ | \ c_i=0, v_i>0, \ \ \forall i \in \scrN,\notag \\ 
    & \hspace{10em} \text{ and} \sum_{i} v_{i}=m-m_{\bu} \}.
\end{flalign}

\begin{definition}[Locally asymptotically stable] \label{def:local-asyp}
    We say that a set of equilibria $\Upsilon_{\bu}$ is locally asymptotically stable if for any equilibrium $(\underline{\bc}, \underline{\bv})\in \Upsilon_{\bu}$, there exists a neighborhood  $\scrB_{\bu}^\delta(\underline{\bc}, \underline{\bv})$  such that every evolution of the model \eqref{eq:system-dynamics} starting at $(\bc(\tau), \bv(\tau))=(\underline{\bc},\underline{\bv})$, and with $(\bc(0), \bv(0))\in \scrB_{\bu}^\delta(\underline{\bc}, \underline{\bv})$ has a limit which belongs to the equilibrium set $\Upsilon_{\bu}$ i.e., $(\lim_{t\xrightarrow{} +\infty} \bc(t), \lim_{t\xrightarrow{} +\infty} \bv(t) )\in \Upsilon_{\bu}$, where $\tau \in [ -\max_{i,j}T_{ij},0)$ and
    \begin{flalign}
        &\scrB_{\bu}^\delta(\underline{\bc}, \underline{\bv}) :=\{ (\bc, \bv) \in \mathbb{R}^{2N} \ | \ c_i>0, v_i=\underline{v_i}, \ \forall i \in \scrN , \notag \\
        & \hspace{10em} \text{ and } ||(\bc-\underline{\bc}, 0)|| < \delta) \}.
    \end{flalign}
\end{definition}

\begin{theorem}[Stability of the equilibria] \label{thm:stability-eq}
    Let $\bu \in \scrU$ and assume $m_{\bu}>\underline{m}$; then, the set of equilibria $\Upsilon_{\bu}$ is locally asymptotically stable.
\end{theorem}
\begin{pf}
    Provided in the Appendix \qed{}
\end{pf}

    \section{Optimal Strategies} \label{sec:optimal-pricing}

In this section, we present an optimization framework to find optimal prices given endogenous demand rates. This model aims to maximize the revenue of an AMoD provider while ensuring load balancing of clients and vehicles. We then turn to a model which uses a rebalancing formulation to ensure load balancing, without the need of price adjustments. Finally, we combine these two formulations into a joint model. 

\subsection{Optimal Pricing} \label{subsec:optimal-pricing}

We are interested in finding the best pricing policy that ensures the existence of an equilibrium \eqref{eq:equilibria}. Hence, we define the feasible set of the pricing problem to be
\begin{flalign*}
    & \scrF = \Big\{\bu  :  \sum_i \big( \lambda_{ij}(u_{ij}) - \lambda_{ji}(u_{ji})\big)=0, \\[-1.2em] 
    & \hspace{14em}\forall j \in \scrN, \ \bu \in [1,\bu^{\max}] \Big\}.
\end{flalign*}
Then, we can define the profit maximization problem as 
\begin{align} \label{eq:pricing-problem}
    &\max\limits_{\bu \in \scrF}  \quad \sum\limits_{ij} \lambda_{ij}(u_{ij})u_{ij}p^0_{ij}-c^o_{ij}\lambda_{ij}(u_{ij})\notag \\[-0.4cm]
    &\hspace{3cm} -c^c(\lambda^0_{ij}(u_{ij}) -  \lambda_{ij}(u_{ij})),
\end{align}
where $\lambda_{ij}(u_{ij})u_{ij}p^0_{ij}$ and $c^o_{ij}$ are the total revenue and the operational cost of request $i$ to $j$, respectively; and $c^c$ is the penalty that the AMoD service incurs when a costumer exists the platform because of a high price. 

Note that if the functions $J_{ij}(u_{ij}):=\lambda_{ij}(u_{ij})u_{ij}p^0_{ij}-c^o_{ij}\lambda_{ij}(u_{ij})-c^c(\lambda^0_{ij}(u_{ij}) -  \lambda_{ij}(u_{ij}))$ are concave in the range of $[1,\bu^{\max}]$, then the optimization problem is tractable (we maximize over a concave function with linear equality constraints). To ensure the concavity of the cost function $J_{ij}$ we  need its second derivative to satisfy
\begin{equation} \label{eq:concavity-condition}
     \ddot{J_{ij}} \leq 0 \implies \ddot{\lambda}_{ij}(u_{ij}) \leq - \frac{2}{u_{ij}p^0_{ij}-c^o_{ij}-c^c} \dot{\lambda}_{ij}(u_{ij}).
\end{equation}
 
 Recall that by Assumption $1$ ($\lambda_{ij}$ is monotonically decreasing) $\dot{\lambda}_{ij}<0$. Hence, for any linear demand function, the problem becomes tractable. 

\subsection{Optimal Rebalancing}
We use the planning rebalancing model developed in ~\cite{pavone2012}. In this setting, we aim to find a static rebalancing policy that reaches an equilibrium. Let the rebalancing flow be $r_{ij}$, that is, the rate at which empty vehicles flow from $i$ to $j$. To solve the problem we use the following Linear Program (LP) that minimizes the empty travel time and seeks to equate the inflow and outflow of vehicles at each region by using $N^2$ variables 
\begin{subequations} \label{eq:rebalancing-problem}
    \begin{flalign}
        \min_{\br \geq 0}& \quad  \sum_{ij} T_{ij}r_{ij} \label{eq:rebalancing-problem-obj} \\
        &\text{s.t. } \sum\limits_{i} \lambda^0_{ij} + r_{ij} -\lambda^0_{ji}- r_{ji} = 0, \hspace{2mm} \forall j \in \scrN  \label{eq:rebalancing-problem-cnstrs}.
    \end{flalign}
\end{subequations}
Notice that in this case we use $\lambda^0_{ij}$ instead of $\lambda_{ij}(u_{ij})$ as we do not consider the possibility of decreasing the demand by using price incentives.
This LP is always feasible as one can always choose $r_{ij}=\lambda^0_{ji}>0$ for all $i,j \in \scrN$ which satisfies the set of constraints \eqref{eq:rebalancing-problem-cnstrs}. 
All the results presented in Section \ref{sec:well-poss-eq-stability} hold for this problem as well and are studied in~\cite{pavone2012}.

\subsection{Joint Pricing and Rebalancing} \label{subsec:joint-pricing-rebalancing}
We are interested in choosing the best policy which leverages different decisions that the mobility-on-demands providers face. In particular, we would like to optimize the pricing, re-balancing and sizing problem. Then, we can write the planning optimization problem as, 

\begin{subequations} \label{eq:problem-full}
    \begin{align} 
    \max\limits_{\bu, \br, m} & \hspace{2mm} \sum\limits_{ij} 
    \lambda_{ij}(u_{ij})u_{ij}p^0_{ij}-c^o_{ij}\lambda_{ij}(u_{ij})  \label{eq:problem-full-obj} \\[-0.4cm] 
     &\hspace{1cm} -c^c(\lambda^0_{ij}(u_{ij}) -\lambda_{ij}(u_{ij}))
    - c^r(r_{ij}T_{ij}) -c^f m \notag \\[0.2cm]
    \text{s.t.} & \hspace{2mm}\sum\limits_{i} \lambda_{ij}(u_{ij})+r_{ij}-\lambda_{ji}(u_{ji})-r_{ji} = 0, \label{eq:problem-full-eq}\\[-.5cm]
    & \hspace{5.8cm}\forall j \in \scrN \notag\\
    & \hspace{2mm} \sum\limits_{ij}
    T_{ij}(\lambda_{ij}(u_{ij})+r_{ij}) \leq m, \\
    & \hspace{2mm} \bu \in [1,\bu^{\max}],
\end{align}
\end{subequations}
where $c^r$ and $c^f$ are the cost of rebalancing and the cot of owning and maintaining a vehicle per unit time, respectively.  Note that to ensure that solving \eqref{eq:problem-full} reaches a global maximum, we must validate that \eqref{eq:concavity-condition} holds for $\bu \in [1,\bu^{\max}]$.

Note that this problem, if solvable, yields a solution with higher profits than the individual formulations of pricing \eqref{eq:pricing-problem} and \eqref{eq:rebalancing-problem}, or the sequential approach of solving first the rebalancing problem and then adjusting prices. This happens given that the problem is \emph{jointly} solving for $\bu$ and $\br$ considering simultaneously the full objective of the profit maximization problem \eqref{eq:problem-full-obj}.
    \section{Experiments} \label{sec:Experiments}
We carry out two case studies to assess the benefits of solving the joint problem of pricing and rebalancing over other approaches. Our first experiment uses a fictitious transportation network to analyze sensitivities with respect to the network size. The second one consists of a data-driven case study using historical data from New York City. We report empirical results of the achievable profit of the AMoD system when solving the problem using the different methodologies presented in Table~\ref{tab:policies}. 

\begin{table}[h]
    \centering
    \caption{Different policies evaluated to plan the operation of an AMoD system.}
    \renewcommand{\arraystretch}{1.2}
    \begin{tabular}{ | l | l l | }
        \hline
        \textbf{Policy} & \textbf{Type} & \textbf{Formulation}\\
        \hline
        $\scrP_{ij}$    &   Individual & \eqref{eq:pricing-problem}  \\
        $\scrR_{ij}$    &    Individual & \eqref{eq:rebalancing-problem}   \\
        $\scrR_{ij}\rightarrow\scrP_{ij}$ &    Sequential &    \eqref{eq:rebalancing-problem} then \eqref{eq:pricing-problem}\\
        $\scrP_{i}+\scrR_{ij}$  & \begin{tabular}{@{}c@{}}Joint with fixed \\[-.5em] price by origin\end{tabular}  &  \begin{tabular}{@{}c@{}} \eqref{eq:problem-full} with $u_{ij}=u_{ik}$\\[-.4em] $\forall i,j,k \in \scrN$ \end{tabular} \\
        $\scrP_{ij}+\scrR_{ij}$ &  Joint &  \eqref{eq:problem-full} \\
        \hline
    \end{tabular}
    \label{tab:policies}
    \vspace{-0.2cm}
\end{table}

We begin with the individual the policies  $\scrP_{ij}$ and $\scrR_{ij}$ to see the equilibrium under a policy or rebalancing strategy. We then turn to a \emph{sequential} approach  $\scrR_{ij}\rightarrow\scrP_{ij}$ to solve the problem. Our motivation for this methodology comes from the fact that many corporations tend to have separate pricing and rebalancing departments, which would result in solving the joint problem sequentially. Note that the sequential policy $\scrP_{ij}\rightarrow\scrR_{ij}$ is not included because once the pricing problem is solved, the system is at equilibrium and the rebalancing problem becomes trivial (i.e., $\br=0$). Finally, the \emph{joint with fixed prices by origin} policy $\scrP_{i}+\scrR_{ij}$ is motivated by the fact that current MoD services only use the origin (not the destination) when setting surge prices (price multipliers)~\cite{chen2015peeking,cohen2016using}.

Note that in order to have a tractable solution for formulations \eqref{eq:pricing-problem} and \eqref{eq:problem-full} we require a function satisfying $\eqref{eq:concavity-condition}$. To achieve this, we assume a linear demand (willingness-to-pay) function, specifically we let
\begin{equation}\label{eq:linear-demand-function}
    \lambda_{ij}(u_{ij})=\frac{\lambda^0_{ij}}{u^{\max}_{ij}-1}(u^{\max}_{ij}-u_{ij}),
\end{equation}
where we set $u^{\max}_{ij}=4$ as suggested in~\cite{cohen2016using}. Hence, by using this linear demand function, we get a tractable Quadratic Program (QP) with linear constraints. Arguably, linear demand functions may not be as accurate as desired for realistic implementations of this model. However, using linear functions allows us to recover a global maximum solution to the problem and assess the potential benefits that \emph{joint} policies may achieve compared to other strategies.

For both experiments we let the \emph{base} price be proportional to the travel time using  $p^0_{ij}=0.5T_{ij}$, where $\$0.5$ is the average price a user pays in dollars per minute of taxi ride reported in~\cite{taxi_fares}. Additionally, we let the operation and rebalancing cost per kilometer be $c^o=c^r=\$0.72$ as suggested in \cite{BOSCH201876};  the lost customer cost $c^c$ is equal to $\$5$,  and we set the reguralizer parameter on the fleet size to be $c^f=\$1\times10^{-10}$.

\subsection{Uniform Demand}
We compare the solution of the different methods for a network with random uniform demands. For each strategy, we let the \emph{base} demand be $\lambda_{ij}^0\sim U(0, 4)$ and travel time between regions be $T_{ij}\sim U(0, 40)$. Then, we solve the problem for networks with a number of regions ranging between $0$ and $60$. 

\begin{figure}[t]
    \centering
    \begin{subfigure}{.41\linewidth}
      \centering
      \includegraphics[trim={0 0cm 0cm 0}, width=1\linewidth]{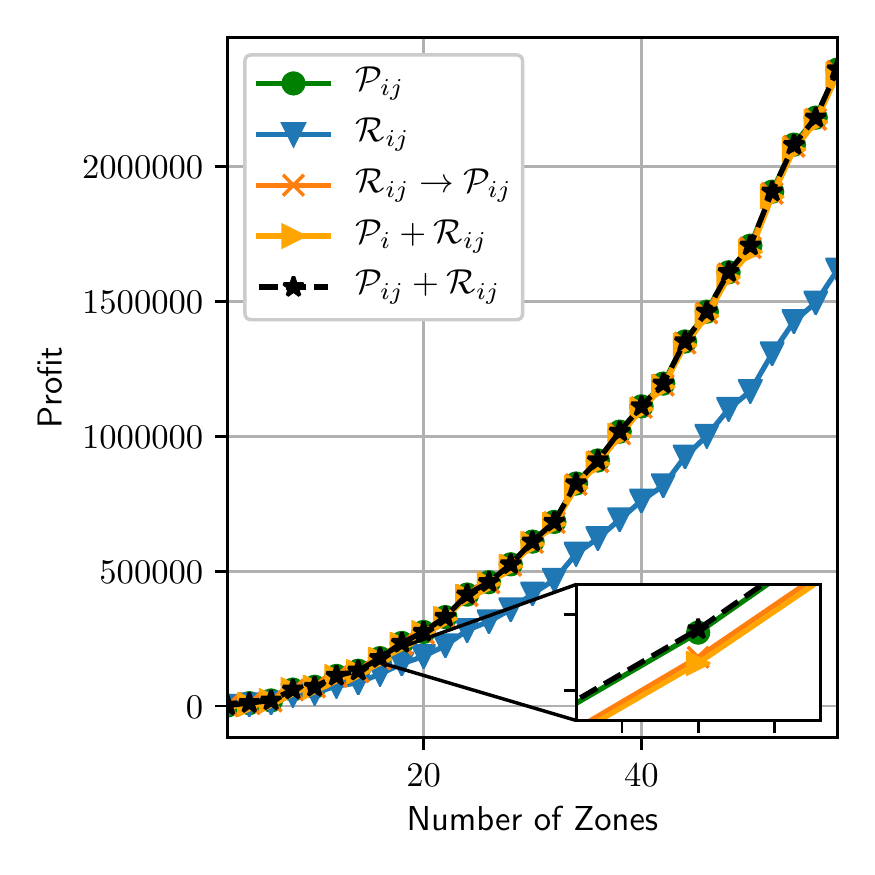}
      \vspace{-0.7cm}
      \caption{}
      \label{fig:price-absolute}
    \end{subfigure}
    \begin{subfigure}{.57\linewidth}
      \centering
      \includegraphics[trim={0cm 0cm 0cm 0}, width=1\linewidth]{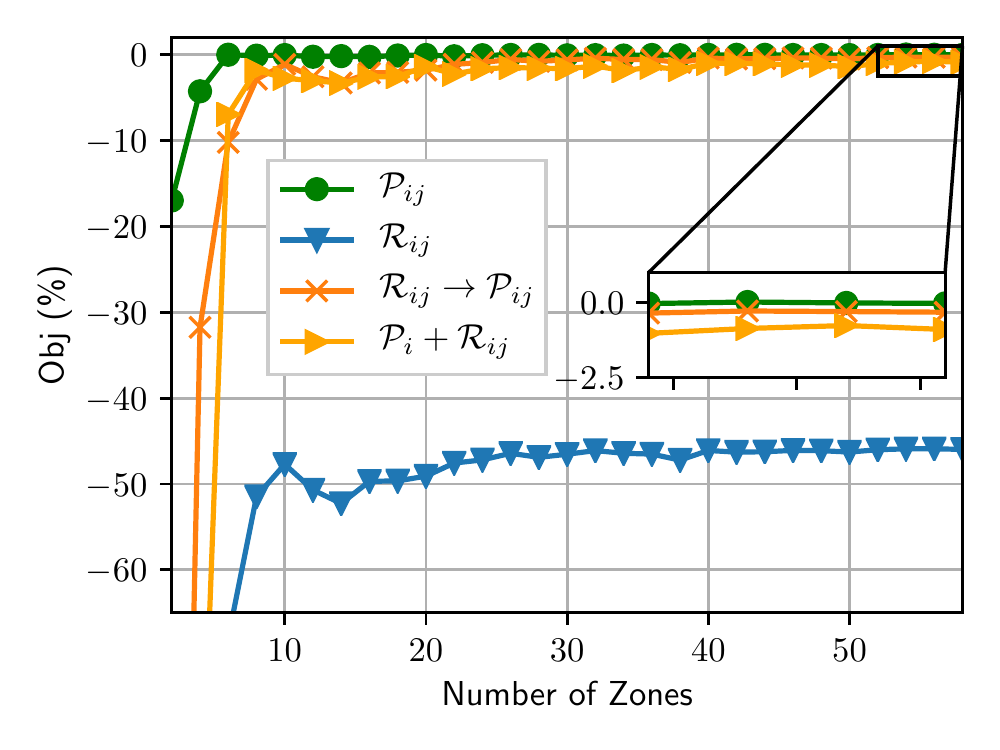}
      \vspace{-0.7cm}
      \caption{}
      \label{fig:price-relative}
    \end{subfigure}
    \caption{Objective function value under different number of zones and AMoD strategies, (a) Shows absolute values of \eqref{eq:problem-full-obj} when different policies are implemented while (b) plots the relative difference between the joint solution and the others. }
    \vspace{-0.4cm}
\end{figure}

Figure~\ref{fig:price-absolute} shows the value of the cost function \eqref{eq:problem-full-obj} for each methodology. Moreover, in Figure~\ref{fig:price-relative} we observe the relative deviation in profits for the solution of each strategy against the joint pricing and rebalancing solution. We see that as the number of regions increases, the deviation converges to a stable value. To explain this phenomenon, we define the \emph{potential} of region $i$ to be the load balance deviation when no pricing or rebalancing policy is applied, namely, $\zeta_{i} = \sum_j \lambda_{ij}^0-\lambda_{ji}^0$. Then, since we draw samples from the same uniform distribution to assign all $\lambda_{ij} \ \forall i,j \in \scrN$, the expected value of $\zeta_i$ is equal to zero for all $i$. Hence, this convergence behavior is simply a direct implication of the law of large numbers. Note that, for the same reason, the individual policy $\scrP_{ij}$ converges to zero.

\subsection{New York City Case Study}
We perform a case study of New York city using the data available at~\cite{nyc-trip-data}. 
Specifically, we analyze the data set of \emph{High Volume For-Hire Vehicle Trip Records} of November 2019~\cite{nyc-trip-data}. In order to analyze stable distributions of trips in the network, we filter the data to consider only working days (Monday to Friday). 
Then, we focus on four time slots: Morning Peak (AM) from 7:00-10:00 hrs, Noon (MD) from 12:00-15:00 hrs, Afternoon Peak (PM) from 17:00-20:00 hrs and Night (NT) from 00:00-3:00 hrs. 
For every time window in November 2019, we collect data on origin-destination pairs and travel times of every trip. Then, we compute the average hourly demand and travel times,  and we use these values to preform our analysis and test the different solutions.

Table~\ref{tab:policy-comparison} shows the deviation in profits (in percentage terms) between the different approaches and the joint formulation. As a reminder, use Table~\ref{tab:policies} summarizes all policy definitions.
We observe that the joint method outperforms all the other methods in the range from $5\%$ to $40\%$, highlighting the benefit of solving this problem using a joint strategy. 
In particular, we observe that each of the individual strategies performs on average worse than strategies that optimize both pricing and rebalancing. 
Also, it is relevant to stress the $5\%$ deviation of the policy with \emph{fixed surge price by origin}, as it matches our expectations of the relevance of considering the destination when pricing.

\begin{figure}[t]
    \centering
    \includegraphics[width=1\linewidth]{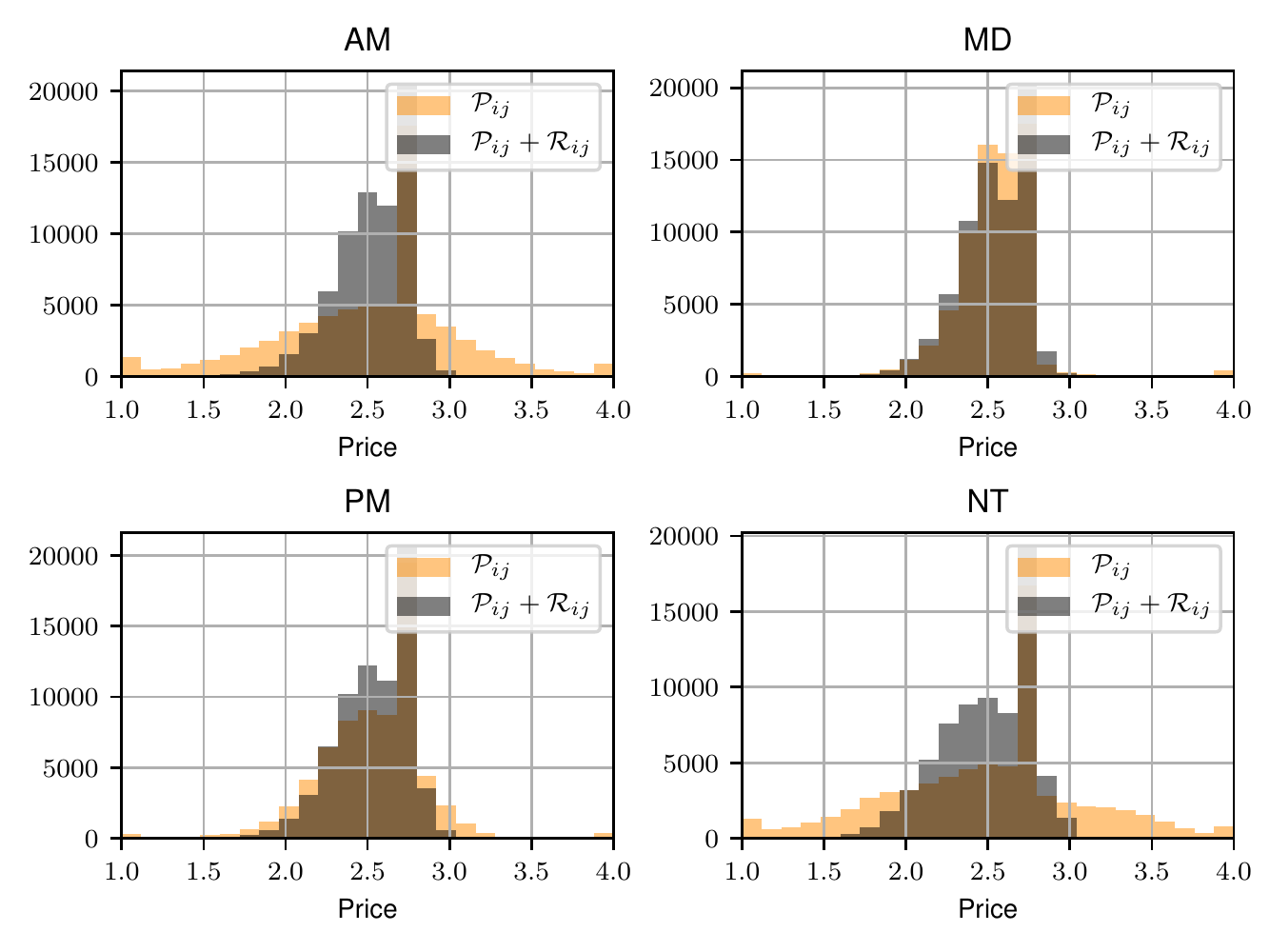}
    \vspace{-0.7cm}
    \caption{Distribution of prices $\bu^*$ for different policies at different time slots}
    \vspace{-0.2cm}
    \label{fig:histogram-prices}
\end{figure}

\begin{table}[t]
    \centering
    \caption{Relative deviation in percentage of each policy compared to the joint strategy $\scrP_{ij}+\scrR_{ij}$ for different time slots}
    \begin{tabular}{|l|llll|}
\hline
\textbf{Policy} & \textbf{AM} & \textbf{MD} & \textbf{PM} & \textbf{NT} \\
\hline 

\hline
$\mathcal{P}_{ij}$                            &   -29.83 &    -8.77 &    -6.64 &   -26.00 \\
$\mathcal{R}_{ij}$                            &   -33.33 &   -28.74 &   -29.20 &   -40.67 \\
$\mathcal{R}_{ij}\rightarrow\mathcal{P}_{ij}$ &   -13.72 &    -9.38 &   -10.89 &   -15.75 \\
$\mathcal{P}_{i}+\mathcal{R}_{ij}$            &       -5.3 &       -5.3 &       -5.1 &         -7.0 \\
\hline

\hline
\end{tabular}
    \label{tab:policy-comparison}
    \vspace{-0.4cm}
\end{table}

To better understand the different approaches, we generated plots of the pricing distribution and trend. Figure~\ref{fig:histogram-prices} shows histograms comparing the value of the solution $\bu$ for the individual pricing policy and the joint strategy. 
As expected, we observe the distribution of the individual approach to have higher variance than the joint. This happens given the hard constraint to reach an equilibrium. When no rebalancing is considered as in $\scrP_{ij}$ the policy chooses prices to ensure $\bu \in \scrF$. 
In contrast, when solving the joint problem, the solution leverages rebalancing and pricing and gives the pricing decision more flexibility to concentrate to pick values that maximize profits.

\begin{figure}[t]
    \centering
    \includegraphics[width=1\linewidth]{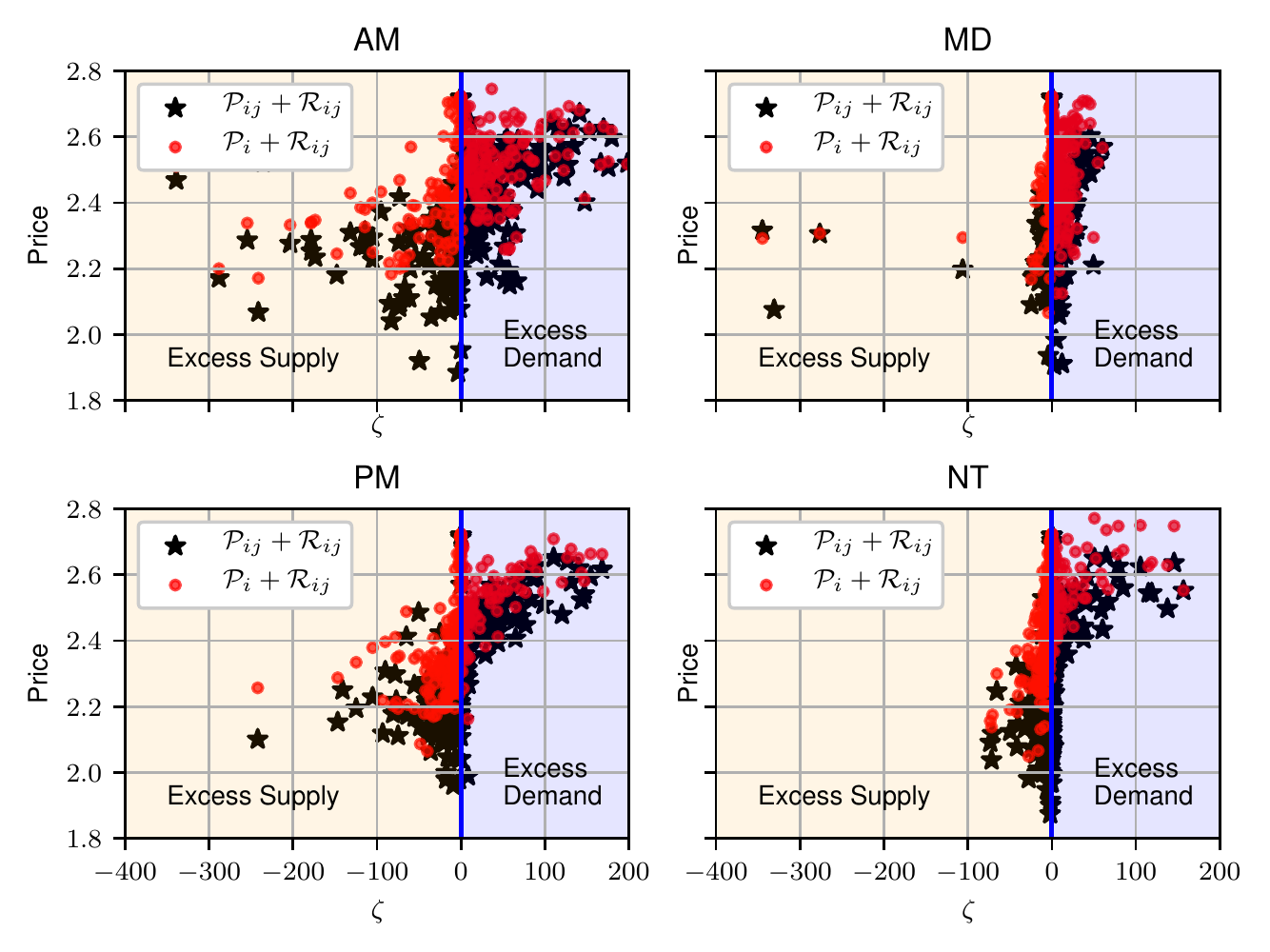}
    \vspace{-0.7cm}
    \caption{$\zeta_i$ indicates excess demand for positive values and excess supply for negative values. For the $\scrP_{i}+\scrR_{ij}$ case we plot prices by origin. For the joint problem $\scrP_{ij}+\scrR_{ij}$, we plot the demand-weighted average price per origin.}
    \vspace{-0.1cm}
    \label{fig:prices-policy}
\end{figure}

Figure~\ref{fig:prices-policy} plots prices against $\zeta_i$ (the  \emph{potential} of origin region $i$) . Recall that parameter $\zeta_i$ indicates excess demand for positive values (i.e., more costumers than vehicles) and excess supply for negative values (i.e., more vehicles than costumers). For the problem with unique surge prices we plot values by origin. For the joint problem we plot the demand-weighted average price by origins. Just as we would expect, there is a positive trend between these variables. The algorithm lowers prices when there is excess supply to incentivize users to request rides, and increases the price when there is excess demand. Note that the pattern is stronger for busier times (AM and PM). 
Finally, we quantify how relevant is the pricing relative to the rebalancing component when balancing the load of the system. Letting $\br^*$ and $\bu^*$ be the solution of \eqref{eq:problem-full}, we define a load dispersion metric as follows $\bar{\zeta}_0 = \frac{1}{N}\sum_i |(\sum_j \lambda_{ij}^0-\lambda_{ji}^0)|$ when nothing is applied, $\bar{\zeta}_\br = \frac{1}{N}\sum_i |(\sum_j \lambda_{ij}^0+r_{ij}-\lambda_{ji}^0-r_{ij})|$ when the rebalancing component is applied, and $\bar{\zeta}_\bu = \frac{1}{N}\sum_i |(\sum_j \lambda_{ij}(u_{ij})-\lambda_{ji}(u_{ji}))|$ when the pricing component (but no rebalancing) is applied. Note that we do not define $\bar{\zeta}_{\bu,\br}$ as the result will be zero given that the system is at equilibrium by \eqref{eq:problem-full-eq}.
Table~\ref{tab:optimal-mean-variance} shows this dispersion metric for the different time slots considered. Interestingly, we see that the pricing component of the policy reduces this metric in all cases, showing its importance for load balancing the system.
\begin{table}[t]
    \centering
    \caption{Dispersion on the average absolute value of potentials when components of the joint policy $\bu^*$ and $\br^*$ are applied.}
    \begin{tabular}{|l|llll|}
\hline

 & \textbf{AM} & \textbf{MD} & \textbf{PM} & \textbf{NT} \\
\hline

\hline

$\bar{\zeta}_0$                            &   57.03 &   16.62 &   34.77 &    17.64 \\

$\bar{\zeta}_{\bu^*}$                     &   20.44 &   4.10 &   6.80 &   6.24 \\

$\bar{\zeta}_{\br^*}$                      &     36.71 &     13.23 &     28.36 &     11.49 \\

\hline

\hline

\end{tabular}

    \label{tab:optimal-mean-variance}
    \vspace{-0.4cm}
\end{table}


    \section{Conclusion}\label{sec:conclusion}
In this paper we studied how a pricing policy which considers origin-destinations can stabilize the system and reach an equilibrium in terms of balancing the load of customer and vehicles. In addition, we formulate a profit maximization optimization model which considers selecting pricing and rebalancing policies jointly. Moreover, we quantify the achievable benefits of solving the problem jointly compared to other methodologies using a data-driven case study of the New York City transportation network. Our results suggest that solving the problem jointly increases the profits of the AMoD provider by up to 40\% when comparing it to individual strategies, 15\% when comparing it to sequential strategies, and 7\% when comparing it to a policy that restricts to a unique \emph{surge} price per origin. 

\emph{Future Work: }
This work can be extended as follows. First, we would like to provide a framework capable of handling more realistic nonlinear demand functions. Second, we would like to complement this model with real-time strategies by the use of a stochastic fluid model~\cite{sun2004perturbation}, as well as a discrete event system~\cite{cassandrasbook} with the aim to provide stochastic and microscopic results of the joint policy. Third, we are interested in coupling this joint solution with the routing problem in \cite{wollensteinbetech2020congestionaware} in order to give an overall optimization framework to operate AMoD systems. Finally, we would like to solve the problem from a welfare maximization perspective rather than from the profit maximization and compare the results. 

    \bibliographystyle{IEEEtran}
    \bibliography{references}

\begin{thebibliography}{10}
\providecommand{\url}[1]{#1}
\csname url@samestyle\endcsname
\providecommand{\newblock}{\relax}
\providecommand{\bibinfo}[2]{#2}
\providecommand{\BIBentrySTDinterwordspacing}{\spaceskip=0pt\relax}
\providecommand{\BIBentryALTinterwordstretchfactor}{4}
\providecommand{\BIBentryALTinterwordspacing}{\spaceskip=\fontdimen2\font plus
\BIBentryALTinterwordstretchfactor\fontdimen3\font minus
  \fontdimen4\font\relax}
\providecommand{\BIBforeignlanguage}[2]{{%
\expandafter\ifx\csname l@#1\endcsname\relax
\typeout{** WARNING: IEEEtran.bst: No hyphenation pattern has been}%
\typeout{** loaded for the language `#1'. Using the pattern for}%
\typeout{** the default language instead.}%
\else
\language=\csname l@#1\endcsname
\fi
#2}}
\providecommand{\BIBdecl}{\relax}
\BIBdecl

\bibitem{wollensteinbetech2020congestionaware}
\BIBentryALTinterwordspacing
S.~Wollenstein-Betech, A.~Houshmand, M.~Salazar, M.~Pavone, C.~G. Cassandras,
  and I.~C. Paschalidis, ``Congestion-aware routing and rebalancing of
  autonomous mobility-on-demand systems in mixed traffic,'' 2020,
  \textit{Submitted to IEEE ITSC Conference 2020}. [Online]. Available:
  \url{https://arxiv.org/pdf/2003.04335.pdf}
\BIBentrySTDinterwordspacing

\bibitem{banerjee2015pricing}
S.~Banerjee, R.~Johari, and C.~Riquelme, ``Pricing in ride-sharing platforms: A
  queueing-theoretic approach,'' in \emph{Proceedings of the Sixteenth ACM
  Conference on Economics and Computation}, 2015, pp. 639--639.

\bibitem{turan2019dynamic}
B.~Turan, R.~Pedarsani, and M.~Alizadeh, ``Dynamic pricing and management for
  electric autonomous mobility on demand systems using reinforcement
  learning,'' 2019.

\bibitem{bimpikis2019spatial}
K.~Bimpikis, O.~Candogan, and D.~Saban, ``Spatial pricing in ride-sharing
  networks,'' \emph{Operations Research}, vol.~67, no.~3, pp. 744--769, 2019.

\bibitem{swaszek2019load}
R.~M. Swaszek and C.~G. Cassandras, ``Load balancing in mobility-on-demand
  systems: Reallocation via parametric control using concurrent estimation,''
  in \emph{2019 IEEE Intelligent Transportation Systems Conference
  (ITSC)}.\hskip 1em plus 0.5em minus 0.4em\relax IEEE, 2019, pp. 2148--2153.

\bibitem{HoerlRuchEtAl2018}
S.~H\"orl, C.~Ruch, F.~Becker, E.~Frazzoli, and K.~W. Axhausen, ``Fleet control
  algorithms for automated mobility: A simulation assessment for {Z}urich,''
  2018.

\bibitem{LevinKockelmanEtAl2017}
M.~W. Levin, K.~M. Kockelman, S.~D. Boyles, and T.~Li, ``A general framework
  for modeling shared autonomous vehicles with dynamic network-loading and
  dynamic ride-sharing application,'' \emph{Computers, Environment and Urban
  Systems}, vol.~64, pp. 373 -- 383, 2017.

\bibitem{ZhangPavone2016}
R.~Zhang and M.~Pavone, ``Control of robotic {Mobility-on-Demand} systems: A
  queueing-theoretical perspective,'' vol.~35, no. 1--3, pp. 186--203, 2016.

\bibitem{IglesiasRossiEtAl2016}
R.~Iglesias, F.~Rossi, R.~Zhang, and M.~Pavone, ``A {BCMP} network approach to
  modeling and controlling {Autonomous} {Mobility-on-Demand} systems,'' 2016.

\bibitem{PavoneSmithEtAl2012}
M.~Pavone, S.~L. Smith, E.~Frazzoli, and D.~Rus, ``Robotic load balancing for
  {Mobility-on-Demand} systems,'' vol.~31, no.~7, pp. 839--854, 2012.

\bibitem{RossiZhangEtAl2017}
F.~Rossi, R.~Zhang, Y.~Hindy, and M.~Pavone, ``Routing autonomous vehicles in
  congested transportation networks: Structural properties and coordination
  algorithms,'' vol.~42, no.~7, pp. 1427--1442, 2018.

\bibitem{SalazarTsaoEtAl2019}
M.~Salazar, M.~Tsao, I.~Aguiar, M.~Schiffer, and M.~Pavone, ``A
  congestion-aware routing scheme for autonomous mobility-on-demand systems,''
  2019.

\bibitem{pavone2012}
M.~Pavone, S.~L. Smith, E.~Frazzoli, and D.~Rus, ``Robotic load balancing for
  mobility-on-demand systems,'' \emph{The International Journal of Robotics
  Research}, vol.~31, no.~7, pp. 839--854, 2012.

\bibitem{filippov2013differential}
A.~F. Filippov, \emph{Differential equations with discontinuous righthand
  sides: control systems}.\hskip 1em plus 0.5em minus 0.4em\relax Springer
  Science \& Business Media, 2013, vol.~18.

\bibitem{haddad1981monotone}
G.~Haddad, ``Monotone trajectories of differential inclusions and functional
  differential inclusions with memory,'' \emph{Israel journal of mathematics},
  vol.~39, no. 1-2, pp. 83--100, 1981.

\bibitem{chen2015peeking}
L.~Chen, A.~Mislove, and C.~Wilson, ``Peeking beneath the hood of uber,'' in
  \emph{Proceedings of the 2015 internet measurement conference}, 2015, pp.
  495--508.

\bibitem{cohen2016using}
P.~Cohen, R.~Hahn, J.~Hall, S.~Levitt, and R.~Metcalfe, ``Using big data to
  estimate consumer surplus: The case of uber,'' National Bureau of Economic
  Research, Tech. Rep., 2016.

\bibitem{taxi_fares}
\BIBentryALTinterwordspacing
``Taxi fares: How much does a ride cost?'' [Online]. Available:
  \url{https://www.seattle.gov/your-rights-as-a-customer/file-a-complaint/taxi-for-hire-and-tnc-complaints/taxi-fares-how-much-does-a-ride-cost}
\BIBentrySTDinterwordspacing

\bibitem{BOSCH201876}
P.~M. Bösch, F.~Becker, H.~Becker, and K.~W. Axhausen, ``Cost-based analysis
  of autonomous mobility services,'' \emph{Transport Policy}, vol.~64, pp. 76
  -- 91, 2018.

\bibitem{nyc-trip-data}
\BIBentryALTinterwordspacing
``Tlc trip record data.'' [Online]. Available:
  \url{https://www1.nyc.gov/site/tlc/about/tlc-trip-record-data.page}
\BIBentrySTDinterwordspacing

\bibitem{sun2004perturbation}
G.~Sun, C.~G. Cassandras, Y.~Wardi, C.~G. Panayiotou, and G.~F. Riley,
  ``Perturbation analysis and optimization of stochastic flow networks,''
  \emph{IEEE Transactions on Automatic Control}, vol.~49, no.~12, pp.
  2143--2159, 2004.

\bibitem{cassandrasbook}
C.~G. Cassandras and S.~Lafortune, \emph{Introduction to Discrete Event
  Systems}, 2nd~ed.\hskip 1em plus 0.5em minus 0.4em\relax Springer, 2010.

\end{thebibliography}
    \newpage
    \section*{Appendix}

\subsection*{Proof of Theorem \ref{thm:stability-eq}}

    We start by showing that $v_i(\tau)>c_i(\tau)$ for  $\tau\in[-\max_{i,j}T_{ij}, t)$, which will serve as a key element to the analysis. 
    To do so, we first assume $v_i(\tau)>0$ for all $i\in\scrN$, and observe the system dynamics \eqref{eq:system-dynamics} at time $t$ are
    \begin{subequations}
        \begin{flalign}
            \dot{c_i}(t) & = (\lambda_{i} - \mu_i) H(c_i), \label{eq:system-dynamics-customer-t}\\
            \dot{v_i}(t) &=  - \lambda_{i} + (\lambda_{i} - \mu_i) H(c_i) +  \sum\limits_{j} (\lambda_{ji} \label{eq:system-dynamics-veh-t} \\[-12pt] 
            & \hspace{3em} - (\lambda_{ji} - \mu_j)  H(c^i_j)) \nonumber, \\
            &=  (\lambda_{i} - \mu_i) H(c_i) - \sum\limits_{j} (\lambda_{ji} - \mu_j)  H(c^i_j), \label{eq:system-dynamics-veh-t-2}\\
            & \geq (\lambda_{i} - \mu_i) H(c_i),\label{eq:system-dynamics-veh-t-3} \\
            & = \dot{c}_i(t), \label{eq:system-dynamics-veh-t-4}
        \end{flalign}        
    \end{subequations}
   where all the $H(v_i)$ in \eqref{eq:system-dynamics} are replaced with $1$ since we assume that $v_i(\tau)>0$. The step \eqref{eq:system-dynamics-veh-t-2} is due to the fact that the system is at equilibrium, i.e.  $\sum\limits_{j} \lambda_{ji} - \lambda_{i}=0,\ \forall i \in \scrN$, and the  step \eqref{eq:system-dynamics-veh-t-3}, is a result of $\mu_i>\lambda_i$ which means that $\lambda_i-\mu_i<0$. 
   Given that $\dot{v}_i(t)\geq\dot{c}_i(t)$ and the fact that at the starting point $\bv > \bc$ (i.e., $\underline{\bv} > \mathbf{0}$), we conclude that $v_i(\tau) > c_i(\tau)$ for $\tau \in [-\max_{i,j} T_{ij}, t)$ and $i \in \scrN$.
   
   Two important consequences of this result are that $c_i$ always reaches $0$ before its corresponding $v_i$, and the vehicle time derivative $\dot{v}_i$ is greater than or equal to $0$ after $c_i$ has reached $0$. This follows by observing that all the terms in \eqref{eq:system-dynamics-veh} are all positive when $c_i=0$.
   
  Now we are in a position to show that $\bv(t)>0$ for $t\geq0$. We do this by combining the second consequence in the previous paragraph with the assumption that the initial state of the system is $(\mathbf{0},\underline{\bv})$ and the fact that $\dot{v}_i(\tau)>\dot{c}_i(\tau)$. Thus, since $\bv(t)>0$ for $t\geq0$, we have that $\dot{c}_i(t)=(\lambda_i-\mu_i)H(c_i)\leq 0$ which will implies that $c_i$ will for sure be $0$ when $t\geq T^{\max}$ where $T^{\max}:=\max_{i} \{c_i(0)/(\mu_i-\lambda_i)\}_{i\in \scrN}$ and which implies that $\lim_{t\xrightarrow{} +\infty} \bc(t)=\mathbf{0}$ since both $\dot{c}_i$ and $c_i$ will be equal to $0$ for all $i\in \scrN$.
  
  To show that $\lim_{t\xrightarrow{} +\infty} \bv(t)=\mathbf{v}$  we use the fact that $c_i=\dot{c}_i = 0$ for $t>0$ and insert this into the vehicle dynamics in \eqref{eq:system-dynamics-veh} obtaining $v_i(t)=\lambda_i(v_i)+\sum_j ( \lambda_{ji} H(v^i_j)-(\lambda_{ji}-\mu_{j})H(c^i_j)H(v^i_j))$. Since, $c_i(t)=0$ we observe that after $T^{\max}+\max_{i,j} T_{ij}$ time units $H(c^i_j)$ will be equal to zero and therefore $\dot{v}(t)=0$ for $t>T^{\max}+\max_{i,j} T_{ij}$.
  Moreover, since $\dot{v}_i(t)=0$ the $\lim_{t\xrightarrow{}\infty} v_i(t)$ exists and can be retrieved using $v_i(t) = v_i(0) + \int_0^t \dot{v}_i(s) ds \geq  v_i(0) + \int_0^t \dot{c}_i(s) ds =  v_i(0) + c_i(t)-c_i(0)$. Given that we show that $v_i(0)>c_i(0)$, we conclude that $\lim_{t\xrightarrow{}\infty} v_i(t)>0$. The property $\lim_{t\xrightarrow{}\infty} v_i(t)>0 > m - m_{\bu}$ holds straight from Proposition \ref{prop:well-posedness}.
  
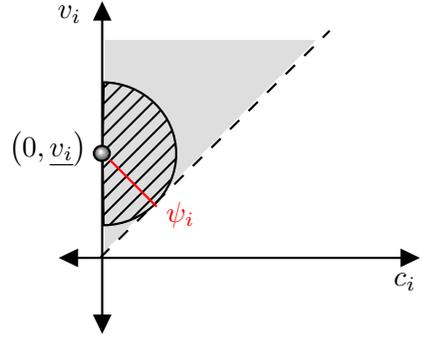
\begin{figure}[t]
    \centering
    \usetikzlibrary{patterns}

 
\tikzset{
pattern size/.store in=\mcSize, 
pattern size = 5pt,
pattern thickness/.store in=\mcThickness, 
pattern thickness = 0.3pt,
pattern radius/.store in=\mcRadius, 
pattern radius = 1pt}
\makeatletter
\pgfutil@ifundefined{pgf@pattern@name@_rvgpdwc4u}{
\pgfdeclarepatternformonly[\mcThickness,\mcSize]{_rvgpdwc4u}
{\pgfqpoint{0pt}{0pt}}
{\pgfpoint{\mcSize+\mcThickness}{\mcSize+\mcThickness}}
{\pgfpoint{\mcSize}{\mcSize}}
{
\pgfsetcolor{\tikz@pattern@color}
\pgfsetlinewidth{\mcThickness}
\pgfpathmoveto{\pgfqpoint{0pt}{0pt}}
\pgfpathlineto{\pgfpoint{\mcSize+\mcThickness}{\mcSize+\mcThickness}}
\pgfusepath{stroke}
}}
\makeatother

  
\tikzset {_elqfgf0kj/.code = {\pgfsetadditionalshadetransform{ \pgftransformshift{\pgfpoint{89.1 bp } { -108.9 bp }  }  \pgftransformscale{1.32 }  }}}
\pgfdeclareradialshading{_whoa6p5u4}{\pgfpoint{-72bp}{88bp}}{rgb(0bp)=(1,1,1);
rgb(0bp)=(1,1,1);
rgb(25bp)=(0,0,0);
rgb(400bp)=(0,0,0)}
\tikzset{every picture/.style={line width=0.75pt}} 

\begin{tikzpicture}[x=0.75pt,y=0.75pt,yscale=-1,xscale=1]

\draw  [color={rgb, 255:red, 255; green, 255; blue, 255 }  ,draw opacity=0.02 ][fill={rgb, 255:red, 222; green, 222; blue, 222 }  ,fill opacity=1 ] (46.5,131.29) -- (145.5,32.29) -- (46.5,32.29) -- cycle ;
\draw    (29,133) -- (189,133) ;
\draw [shift={(192,133)}, rotate = 180] [fill={rgb, 255:red, 0; green, 0; blue, 0 }  ][line width=0.08]  [draw opacity=0] (8.93,-4.29) -- (0,0) -- (8.93,4.29) -- cycle    ;
\draw [shift={(26,133)}, rotate = 0] [fill={rgb, 255:red, 0; green, 0; blue, 0 }  ][line width=0.08]  [draw opacity=0] (8.93,-4.29) -- (0,0) -- (8.93,4.29) -- cycle    ;
\draw    (46,165) -- (46,18) ;
\draw [shift={(46,15)}, rotate = 450] [fill={rgb, 255:red, 0; green, 0; blue, 0 }  ][line width=0.08]  [draw opacity=0] (8.93,-4.29) -- (0,0) -- (8.93,4.29) -- cycle    ;
\draw [shift={(46,168)}, rotate = 270] [fill={rgb, 255:red, 0; green, 0; blue, 0 }  ][line width=0.08]  [draw opacity=0] (8.93,-4.29) -- (0,0) -- (8.93,4.29) -- cycle    ;
\draw  [dash pattern={on 4.5pt off 4.5pt}]  (45,133) -- (150.5,28.57) ;
\draw  [pattern=_rvgpdwc4u,pattern size=6pt,pattern thickness=0.75pt,pattern radius=0pt, pattern color={rgb, 255:red, 0; green, 0; blue, 0}] (46.66,52.29) .. controls (65.11,52.53) and (80,67.15) .. (80,85.14) .. controls (80,103.15) and (65.09,117.78) .. (46.61,118) -- (46.19,85.14) -- cycle ;
\path  [shading=_whoa6p5u4,_elqfgf0kj] (42.19,84.95) .. controls (42.19,83.04) and (43.82,81.49) .. (45.84,81.49) .. controls (47.86,81.49) and (49.5,83.04) .. (49.5,84.95) .. controls (49.5,86.85) and (47.86,88.4) .. (45.84,88.4) .. controls (43.82,88.4) and (42.19,86.85) .. (42.19,84.95) -- cycle ; 
 \draw   (42.19,84.95) .. controls (42.19,83.04) and (43.82,81.49) .. (45.84,81.49) .. controls (47.86,81.49) and (49.5,83.04) .. (49.5,84.95) .. controls (49.5,86.85) and (47.86,88.4) .. (45.84,88.4) .. controls (43.82,88.4) and (42.19,86.85) .. (42.19,84.95) -- cycle ; 

\draw [color={rgb, 255:red, 255; green, 0; blue, 0 }  ,draw opacity=1 ]   (49.5,87.95) -- (71.05,109.5) ;

\draw (31,21) node    {$v_i$};
\draw (185,144) node    {$c_i$};
\draw (82,113) node  [color={rgb, 255:red, 255; green, 0; blue, 0 }  ,opacity=1 ]  {$\psi _{i}$};
\draw (21,83) node    {$\left(0, \underline{v_i}\right)$};

\end{tikzpicture}
    \caption{Sketch of a variable of the initial solution $(\underline{\bc},\underline{\bv})$ along with its neighborhood $\scrB^\delta_\bu$. Shaded in grey is the feasible region ($c_i<v_i$).}
    \label{fig:psi}
\end{figure}

  Finally, to characterize the ball $\scrB^\delta_{\bu}$ we set $\psi_i:=\underline{v_i}\sin(\pi/4)$ and  $\psi^{\min}:=\min_{i} \psi_i$ (see Fig. \ref{fig:psi}). Then, for $\delta=\psi^{\min}$ any path of the system \eqref{eq:system-dynamics} starting at $(\bc(\tau), \bv(\tau))=(\underline{\bc},\underline{\bv})$  for $\tau\in[-\max_{i,j}T_{ij},0)$ and satisfying $(\bc(0), \bv(0))\in \scrB_{\bu}^\delta(\underline{\bc}, \underline{\bv})$
  has a limit which belongs to the equilibrium set $\Upsilon_{\bu}$.
  \qed{}

\end{document}